\documentclass[preprint,12pt]{elsarticle}
\usepackage{amssymb,graphics,graphicx}
\usepackage{amsmath}

\journal{Mathematics and Computers in Simulation}

\begin{document}

\begin{frontmatter}

\title{A Particle Method for a Collisionless Plasma with Infinite Mass}
\author{Stephen Pankavich}
\address{Department of Mathematics \\
University of Texas at Arlington\\
Arlington, TX 76016 \\
sdp@uta.edu}

\begin{abstract}
The one-dimensional Vlasov-Poisson system is considered and a particle method is developed to approximate solutions without compact support which tend to a fixed background of charge as $\vert x \vert \rightarrow \infty$.  Such a system of equations can be used to model kinetic phenomena occurring in plasma physics, such as the solar wind.  The particle method is constructed, implemented, and used to determine information regarding the time asymptotics of the electrostatic field.
\end{abstract}


\end{frontmatter}

\section*{Introduction}
\noindent The motion of a collisionless plasma - an ionized gas of high-temperature and low-density - is described by the Vlasov-Maxwell system.  If this fundamental kinetic model is posed in a two-dimensional phase space (one for spatial variables and another representing momentum) then Maxwell's equations simplify greatly and the problem reduces to the one-dimensional Vlasov-Poisson system.  We consider this system of equations for the motion of negative charges upon prescribing a fixed, spatially-homogeneous background of positive ions:
\begin{equation}
\label{VP} \left. \begin{gathered}
\partial_t f + v \ \partial_x f - E \ \partial_v f = 0,\\
\partial_x E(t,x) = \int ( F(v) - f(t,x,v) ) \ dv \\
f(0,x,v) = f_0(x,v). \end{gathered} \right \}
\end{equation}
Here, $t \in [0,T]$ denotes time, $x \in \mathbb{R}$ is
space, $v \in \mathbb{R}$ is momentum, $f(t,x,v)$ represents the density of negative ions, $F(v) \not\equiv 0$ is a given function describing the fixed background of positive charge, and $E(t,x)$ represents the electric field generated by the charges. Further define
$$ \rho(t,x) := \int( F(v) - f(t,x,v)) \ dv $$
to be the density of charge in the system.
Unlike many classical investigations of the Vlasov-Poisson system, we are interested in studying how the negative ions move to balance the fixed positive background.  Hence, we seek solutions $f$ which tend to $F$ as $\vert x \vert \rightarrow \infty$, rather than to zero. The introduction of the fixed background implies that the total positive charge, total negative charge, and total energy are all infinite.  In the present study, we will also assume neutrality
$$\int \int (F(v) - f_0(x,v)) \ dx \ dv = 0.$$
This condition then yields zero data at infinity for the electric field, and hence a representation for the field results:
\begin{equation}
\label{E} E(t,x) = \int_{-\infty}^x \rho(t,y) \ dy.
\end{equation}
To set the context of the problem, we impose assumptions on the data.  Specifically, we will assume throughout that $f_0 \in \mathcal{C}^1(\mathbb{R}^2)$ is nonnegative with compact support in $v$, and there is
$R > 0$ such that for $\vert x \vert > R$,
\begin{equation}
\label{IC} f_0(x,v) = F(v)
\end{equation}
where $F \in \mathcal{C}_c^1(\mathbb{R})$ is a nonnegative, even, and decreasing function.

As one can see, no assumption of spatial compact support is made on the initial particle distribution $f_0$.  The impetus for investigating a problem like (\ref{VP}) arises from plasma physicists' attempts to study the stability of a two-species neutral plasma in which a perturbation of the distribution of negative ions is introduced from equilibrium and evolves in such a way as to cancel the effects of the prescribed Maxwellian density of positive charge. Though our assumption of compact velocity support precludes Maxwellian backgrounds $F$, the particle simulation would need to be truncated regardless, and the method will still approximate values of a Maxwellian.  As the initial particle density lacks spatial decay, most of the known mathematical theory regarding existence and uniqueness of solutions does not apply.  Recently, the author \cite{1DLocalExistence} has shown the local-in-time existence and uniqueness of solutions to (\ref{VP}) under similar hypothesis, and this argument can easily be adapted to extend the result to the assumptions above. Having answered the question of local well-posedness, a next logical step is to construct numerical methods that can determine the behavior of solutions.  As such, this is the main objective of the current work.

For problems which use a kinetic description of plasma, particle methods are often utilized to numerically approximate solutions and tend to be much more efficient and accurate than other traditional approximation techniques for partial differential equations, such as finite difference or finite element methods.  While particle methods for the Vlasov-Poisson system have been previously studied [1-6], the introduction of non-zero spatial behavior of the particle density as $\vert x \vert \rightarrow \infty$ eliminates the crucial feature of compact spatial support.  Hence, any truncation of the spatial domain must consider the effects of particles originating from outside this region.  In what follows, we derive and construct a method which combats this problem, while discussing both its abilities and limitations.

The general structure of a particle method can be described quite simply (see \cite{BL}).  As in other numerical methods, phase space (which is $(x,v)$ for kinetic equations) is discretized into grids of finite length. Particles are then initialized with starting positions and
velocities at time $t = 0$.  The charge (and if necessary, current) density is calculated from these particles, and the electric (and if
necessary, magnetic) field is calculated from the density.  Finally, we calculate the force exerted by the field and ``move'' the
particles by changing their respective positions and velocities accordingly.  Since the trajectories have now been calculated for the next time step, the process repeats until a stopping time $t= T$ is reached.  Weighting schemes play a large role in these calculations, specifically because particle charges, positions, and velocities must be recorded ``at the particles'', whereas densities, fields, and forces are indexed by prescribed gridpoints, with the number of particles and gridpoints differing dramatically.

Typically, a particle method tracks the evolution of a finite amount of charge within fixed (or adaptable) spatial and velocity domains. However, for any simulation of (\ref{VP}), which models a plasma density without compact spatial support, the spatial domain must be truncated.  Thus, we choose $L > 0$ and perform the computations on the interval $-L \leq x \leq L$. The assumption on the data is reformulated to lie within the initial spatial domain, i.e. there is $R \in [0,L)$ such that for $R < \vert x \vert \leq L$,
\begin{equation}
\label{DiscreteIC} f_0(x,v) = F(v).
\end{equation}
This ensure that the charges cancel outside of a spatial interval $[-R,R]$. Similarly, the velocity domain must be truncated.  We choose $Q > 0$ so that particles in the simulation may take on
velocities only in the interval $-Q \leq v \leq Q$ at the initial time.  As the process
continues, however, the velocity domain is enlarged to allow the particles to move as dictated
by their interaction with the self-consistent electric field.   The compact velocity support of $f$ and $F$ ensure that the largest particle velocity attained at any timestep is finite.  Hence, at every timestep this value is computed and used to extend the spatial boundary of particle dynamics.  Thus, the spatial domain is also enlarged with time, dependent upon current velocities of the particles, to ensure that none can escape.

In the present context, particles are allowed to move neither into nor out of the computational domain. Though the positive and negative charges cancel outside of the spatial interval $[-R,R]$, both positive and negative particles exist outside of this domain and may influence the computations.  Hence, in addition to enlarging the spatial grid, the domain of validity, in which the particles beginning inside the interval $[-L,L]$ could not have been influenced by those which began outside, is computed after the simulation is complete.  The observed values of $\rho$ and $E$ are only considered valid inside this region.

\section*{Description of the Method}
\noindent In this section the particle method is constructed.  A leap-frog scheme is utilized for the particle trajectories and first-order averaging methods are used to interpolate the field and charge density values. Begin by choosing $\Delta x, \Delta v > 0$ and define for every $i,j \in \mathbb{Z}$,

$$ X_{ij}(0) = i \cdot \Delta x $$
$$ V_{ij}(0) = j \cdot \Delta v $$
$$ q_{ij} = f_0(X_{ij}(0), V_{ij}(0)) \cdot \Delta x \ \Delta v.$$
These quantities represent the initial particle positions, initial particle velocities, and the total negative charge included in the simulation, respectively. The functions  $X_{ij}(t)$ and $V_{ij}(t)$ will be defined later for $t > 0$.  Once they are known the approximation of the continuous number density is then given by

\begin{equation}
f(t,x,v) = \sum_{i,j} q_{ij} \hat{\delta}(x - X_{ij}(t)) \ \delta(v - V_{ij}(t))
\end{equation}
where $\delta$ is the Dirac mass and $\hat{\delta}$ is the first-order weighting function defined by
$$ \hat{\delta}(x) = \left \{ \begin{array}{rl} \displaystyle \frac{1}{\Delta x}\left (1 - \frac{\vert x \vert}{\Delta x} \right ),  & \mathrm{if} \ \vert x \vert < \Delta x \\  0, & \mathrm{otherwise.} \end{array} \right.$$
Choose $\Delta t > 0$ and define $t^n = n \cdot \Delta t$ and $x_l = l \cdot \Delta x$ for $n
\in \mathbb{N}$, $l \in \mathbb{Z}$.  We will write $$ E_l^n = E(t^n, x_l)$$ for field values at the $n^{th}$ timestep and $l^{th}$ gridpoint and define the function $E^n(x)$ by
linear interpolation of the gridpoint values $E_l^n$.

To initiate the leap-frog scheme, we first shift the particle velocities backward by a half timestep using the initial field values.  Hence, let $$V_{ij}(t^{-1/2}) = V_{ij}(0) - E(X_{ij}(0)) \cdot \frac{\Delta t}{2}.$$ Now, for $n \in \mathbb{N} \cup \{0\}$, assuming $V_{ij}(t^{n-1/2})$ has been computed, define $$ S^{n-1/2} := \sup_{i,j \in \mathbb{Z}} \left \vert V_{ij}(t^{n-1/2}) \right \vert$$
to be the largest velocity at time $t^{n -1/2}$.  The length
of the spatial domain at time zero, denoted by $L$, must be enlarged at every
time step, so define $L^0 = L$ and for every $n \in \mathbb{N}$,
$$ L^n := L^{n-1} + S^{n-1/2} \cdot \Delta t.$$
Thus, $L_n$ will be the length of the spatial grid at the $n^{th}$
time step.  This will enable the program to continually account for
particles which began inside the inital spatial domain, but due to
an increase in their velocity would normally move outside of this
region.  By enlarging the spatial domain at every time step in this
fashion, we are able to track such particles and their resulting effects on the induced
field.

\input{epsf}
\begin{figure}[t]
\epsfysize=3in
\centering \includegraphics[scale=1.4]{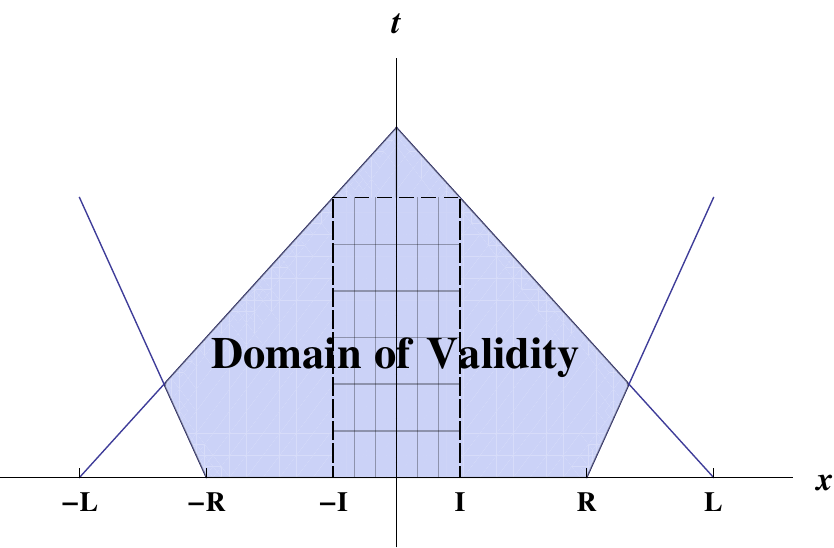}
\caption{The domain of
validity (shaded) and the chosen long-time region (dashed) for the particle method approximations} \label{DOV}
\end{figure}

Assume for $n \in \mathbb{N}$, $X_{ij}(t^n)$ and $V_{ij}(t^{n-1/2})$ are known for all $i,j \in
\mathbb{Z}$, and that $E^n_l$ is known for $\displaystyle \vert l \vert \leq \frac{L_n}{\Delta x}$.  For $ \displaystyle \vert l
\vert > \frac{L_n}{\Delta x}$, in light of the neutrality of the plasma, we find $ E^n_l = 0$.  After a linear interpolation of the field values at gridpoints, define $V_{ij}(t^{n+1/2})$ by
$$ \frac{V_{ij}(t^{n+1/2}) - V_{ij}(t^{n - 1/2})}{\Delta t} = q_{ij} \cdot E^n (X_{ij}(t^n) ).$$
and $$ X_{ij}(t^{n+1}) = X_{ij}(t^n) + \Delta t \cdot V_{ij}(t^{n + 1/2}).$$
Next, we describe how to advance the electric field.  Define for $l
\in \mathbb{Z}$ and $n \in \mathbb{N}$,
$$ \rho_l^{n+1} = \int \left ( F(v) - f(t^{n+1}, x_l, v) \right ) \ dv,$$
and by linear interpolation of $\rho_l^{n+1}$, define $p^{n+1}(x)$. Then, define
$$ E_l^{n+1} = \int_{-L}^{x_l} \rho^{n+1}(y) \ dy.$$
The process may then continue by advancing the particle positions and velocities, $X_{ij}$ and
$V_{ij}$, to the next time step, $t^{n+1}$.

Finally, as the process continues, we must determine the region on
which our approximations are valid.  Since it was necessary to
truncate the spatial domain, we were forced to neglect the effects
of particles which begin outside of our domain and move with large
enough velocity to enter it at some time step. Therefore, we must
discount the approximations within the region in which these particles
could have entered and affected the computations.  Define, for every
$n \in \mathbb{N}$,
$$ P^n := \sum_{k=0}^n S^{k+1/2}  $$
to be the total sum of largest particle velocities up to time $t^{n+1/2}$.
Since no particle beginning outside the original spatial domain
could move with velocity greater than $P^n$, we may conclude that
the portion of the original spatial domain that is unaffected by
such particles at any time step $t^n$ lies inside the interval $[-L
+ P^n \Delta t, L - P^n \Delta t]$.  As Figure \ref{DOV} illustrates, this
provides us with a specific space-time domain on which our numerical
approximations are valid.

\section*{Validation and Steady States}
\noindent In addition to the situation in which $f_0(x,v) = F(v)$ and hence $E \equiv 0$, the previously described particle method was tested using a known steady state solution defined as follows. Let $$ f(x,v) = \mathcal{F} \left (\frac{1}{2}\vert v
\vert^2 + U(x) \right),$$ where $$ \mathcal{F}(e) = \left \{
\begin{array}{cc} -e & \mathrm{if} \ \ e \leq 0 \\ 0 & \mathrm{if} \
\ e \geq 0 \end{array} \right. $$ and
$$ U(x) := -\frac{1}{2}(1 - x^2)^3 \chi_{(-1,1)}(x).$$
From the steady potential $U(x)$, the resulting time-independent field is calculated by
$$ \mathcal{E}(x) := U'(x) = 3x(1 - x^2)^2 \chi_{(-1,1)}(x).$$
Additionally, the charge density $$\rho(x) = \mathcal{E}'(x) = 3(1-x^2)(1-5x^2) \chi{(-1,1)}(x).$$
Therefore, the distribution of positive charge is determined by $\rho$ and $f$ as
$$F(x,v) = \left ( 3(1 - x^2) (1 - 5x^2)
+ \frac{2}{3}(1 - x^2)^\frac{9}{2} \right ) \chi_{(-1,1)}(x) \delta(v).$$

\begin{table}[t]
\centering
\begin{tabular}{|c|c|c|c|c|c|}
\hline
mesh values \ $\backslash$ \ time values & t = 0 & t = 0.12 & t = 0.24 & t = 0.36 & t = 0.48 \\
\hline
$\Delta t = \Delta x = \Delta v = 0.04$ & $8.0 \times 10^{-3}$ &  $8.0 \times 10^{-3}$ & $0.012$ & $0.016$ & $0.018$ \\
& & & & & \\
$\Delta t = \Delta x = \Delta v$ = $0.02$ & $2.0 \times 10^{-3}$ & $2.0 \times 10^{-3}$ & $3.0 \times 10^{-3}$  & $0.004$ &  $0.0005$ \\
& & & & & \\
$\Delta t = \Delta x = \Delta v = 0.01$ & $5 \times 10^{-4}$ &  $5 \times 10^{-4}$ &  $8 \times 10^{-4}$ & $0.001$ & $0.002$ \\
\hline
$\sup_l \vert E(t,x_l) \vert$ & $0.8548$ & $0.8568$ & $0.8598$ & $0.8620$ & $0.8629$ \\
\hline
\end{tabular}
\caption[Error of test field for 1-D particle method]{Error of the field for steady state solution} \label{TABLE1}
\end{table}

Using these functions, the method was implemented for several choices of $\Delta x,\Delta v$, and $\Delta t$ in order to assure convergence to the correct steady state solution.  Table \ref{TABLE1} summarizes the results of these runs, listing the error found by calculating the difference between the known steady field solution and the computed electric field at every time step.

We expect that as the mesh is refined (and the values of $\Delta t$, $\Delta x$, and $\Delta v$ decrease), the values of the error should also decrease at a suitable rate for each time.  Choosing a time $t$ in the table and evaluating each of the three error values, we see that this is the case, and that as the values of spacings are halved, the error decreases by a factor of $4$.  Thus, we trust the method converges at a quadratic rate and is second order accurate.

\section*{Simulation and Time Asymptotics}
Now, consider the following choices for $F$ and $f_0$. Let $U:\mathbb{R}^2 \times [0,\infty) \rightarrow \mathbb{R}$ be defined by $$ U(z, A, B) := A(B - z^2)^3 \chi_{(-\sqrt{B}, \sqrt{B})}(z),$$ and then define

$$F(v) = U(v,1,1) $$
and
$$ f_0(x,v) = F(v) + x \cdot U(x,1,1) \cdot U(v, 0.1, 0.6).$$

\begin{figure}[t]
\begin{minipage}[h]{0.20\linewidth}
\centering \includegraphics[scale=0.47]{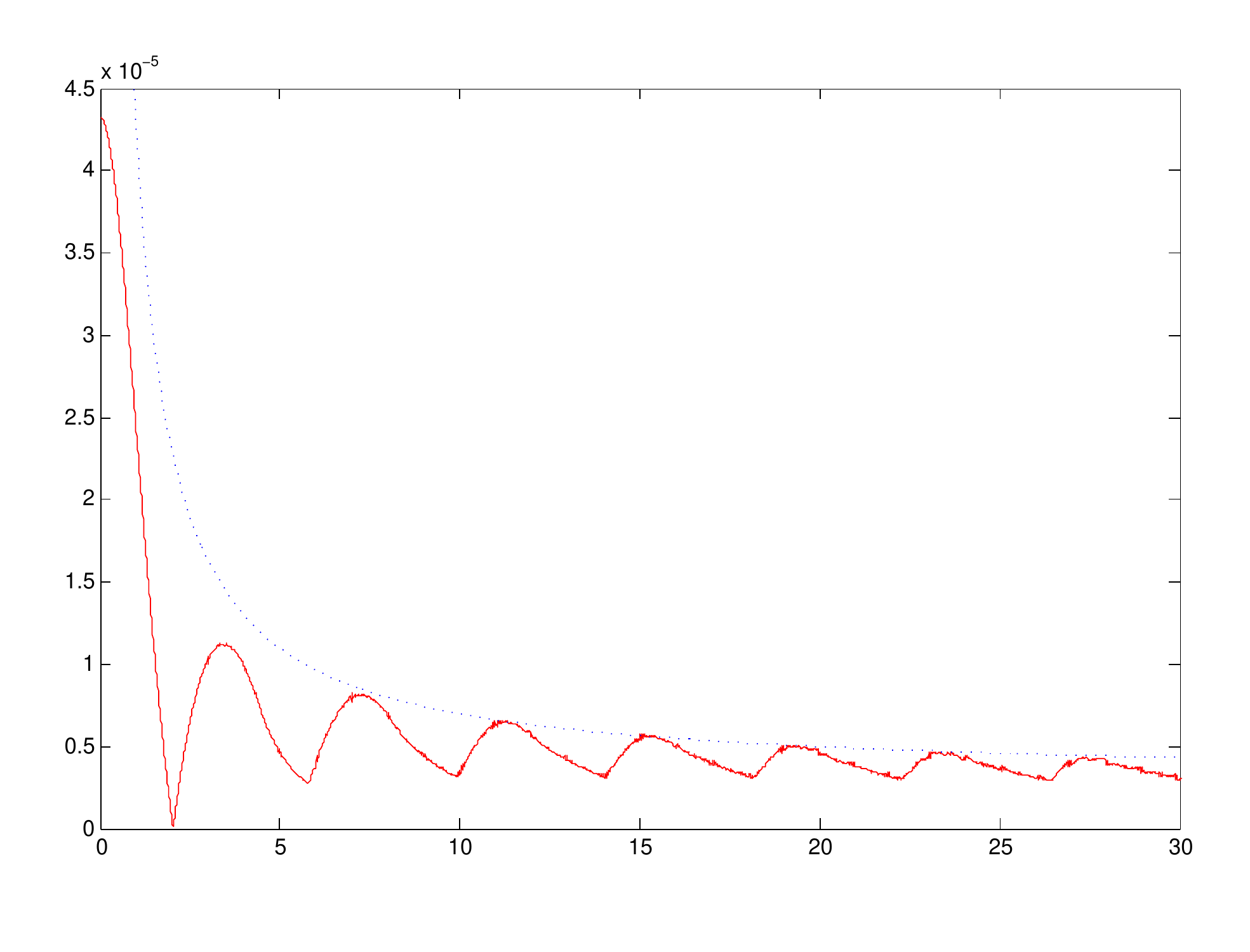}
\end{minipage}\hfill
\begin{minipage}[h]{0.37\linewidth}
\epsfysize=2in
\centering \includegraphics[scale=0.35]{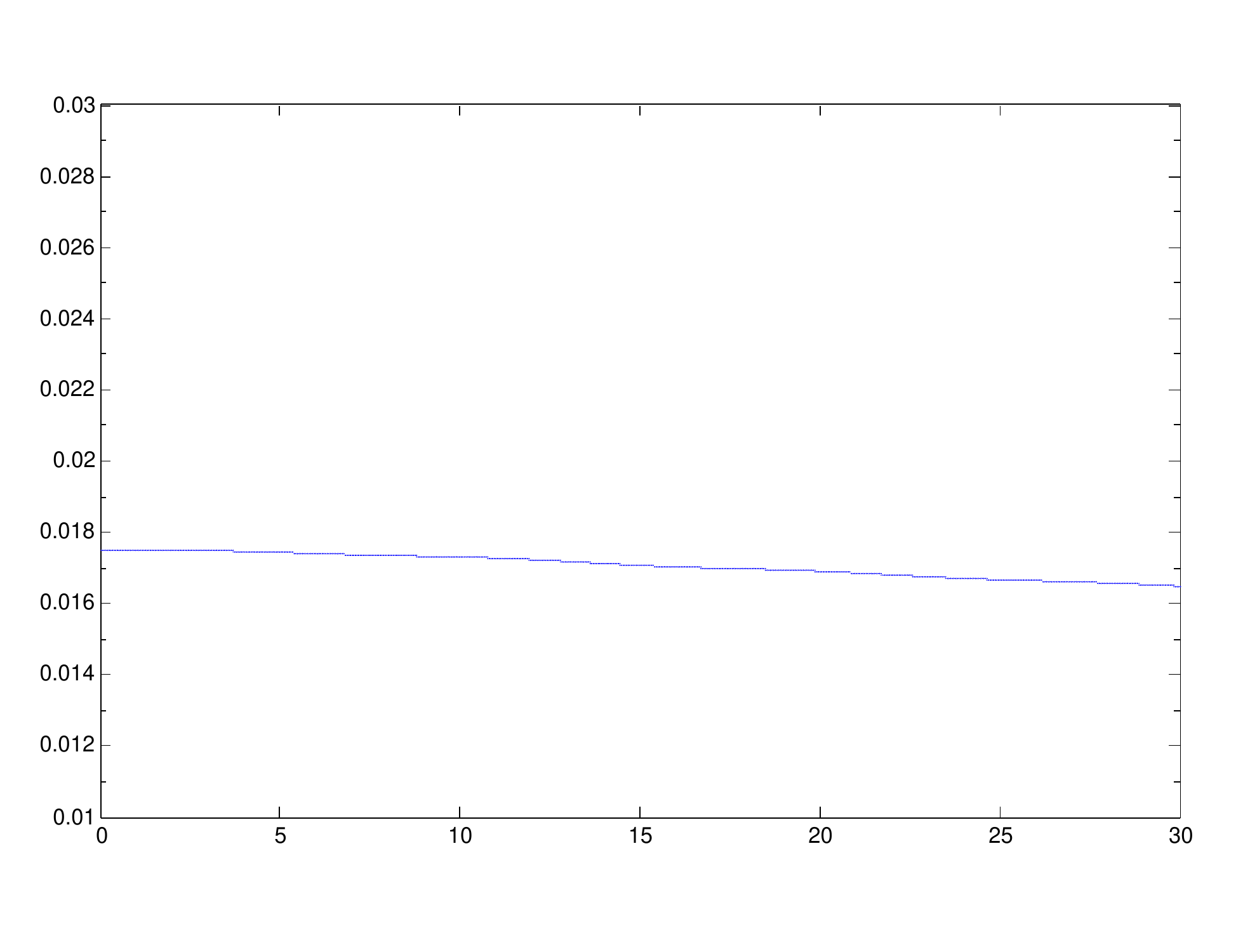}
\end{minipage}
\caption{The computed electric field (left, lined), asymptotic behavior (left, dashed), and net energy (right) for $0 \leq t \leq 30$}
\label{OLDFIELD50}
\end{figure}

\begin{figure}[t]
\begin{minipage}[h]{0.3\linewidth}
\centering \includegraphics[scale=0.65]{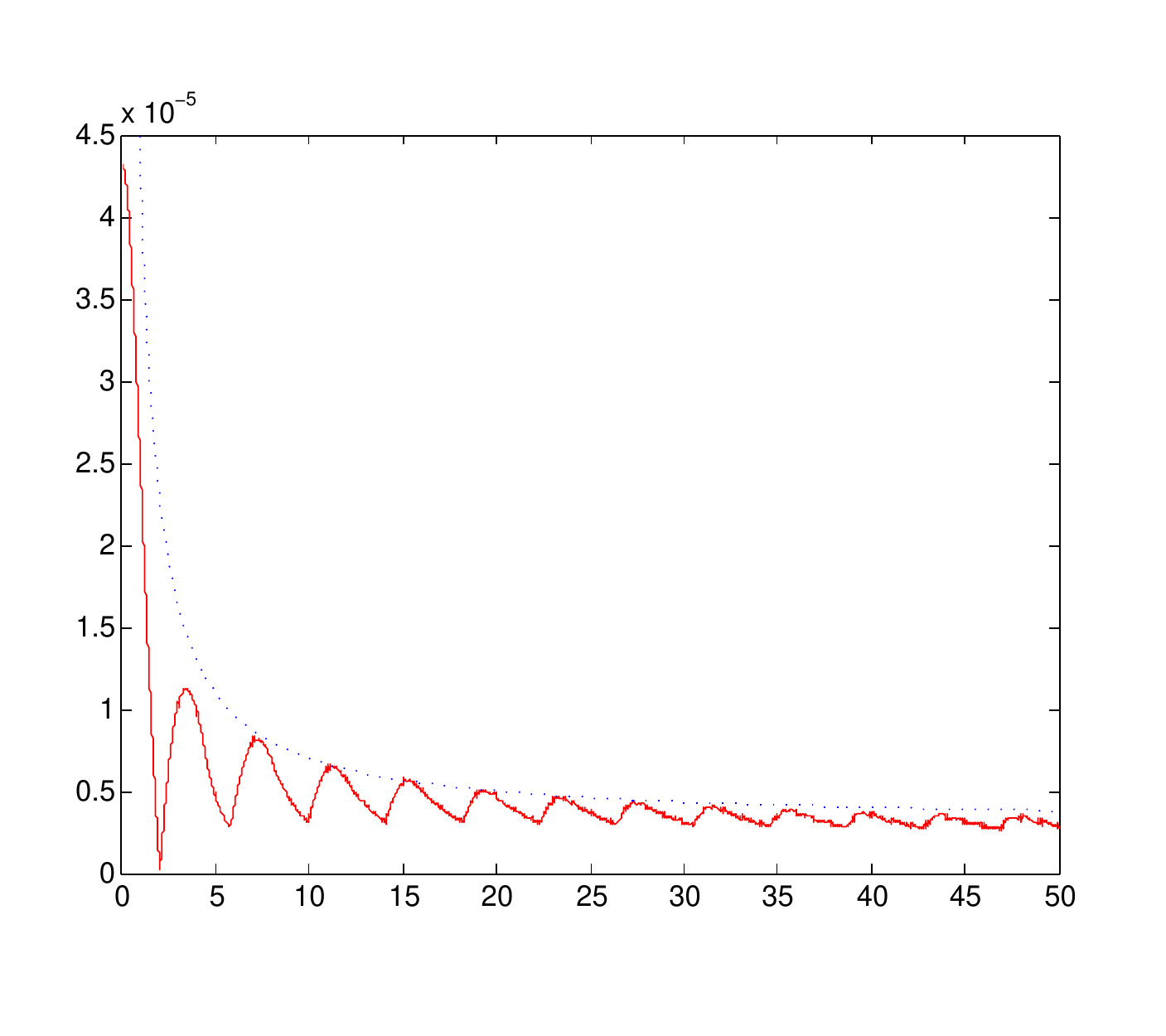} 
\end{minipage}\hfill
\begin{minipage}[h]{0.37\linewidth}
\centering \includegraphics[scale=0.35]{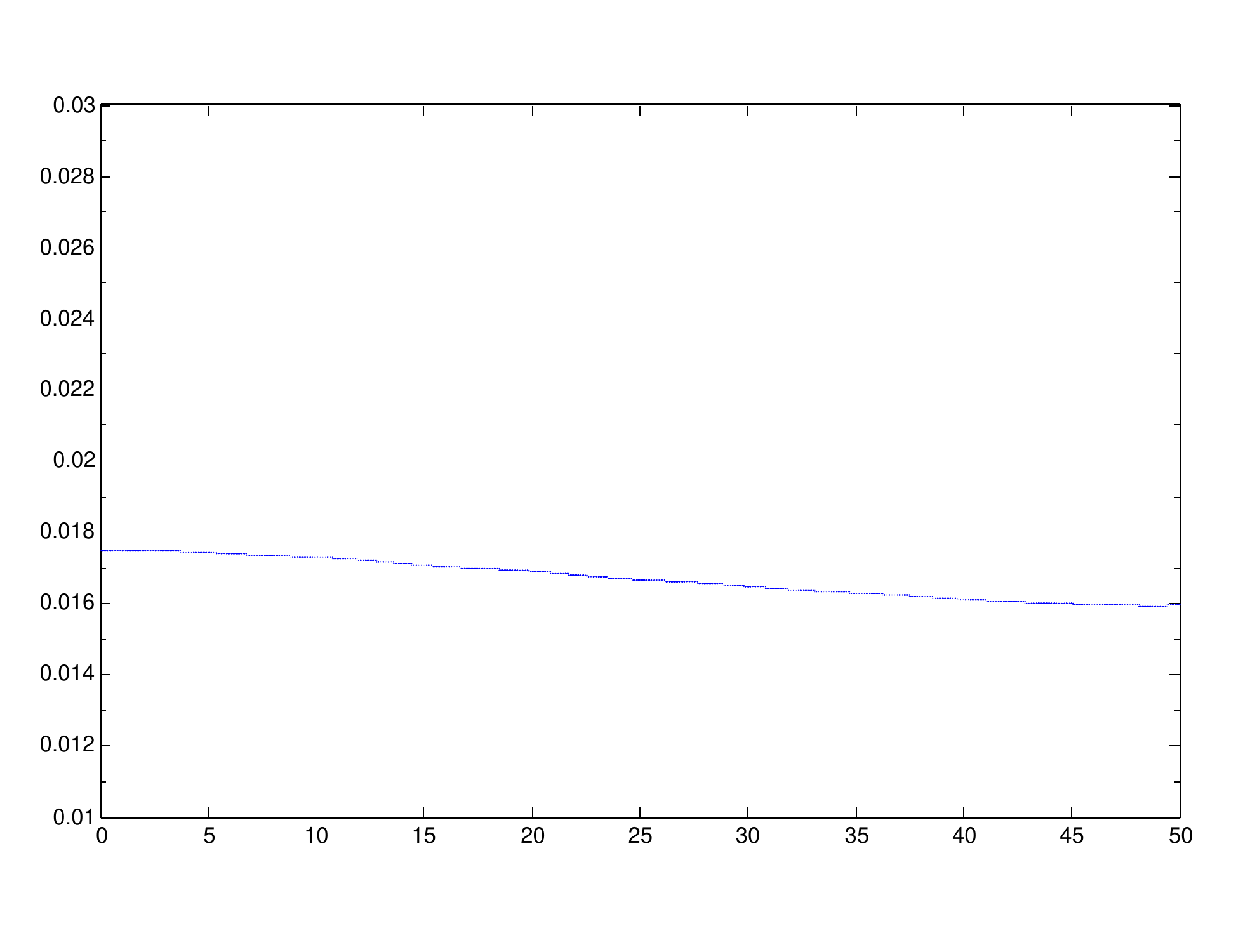} 
\end{minipage}
\caption{The computed electric field (left, lined), asymptotic behavior (left, dashed), and net energy (right) for $0 \leq t \leq 50$}
\label{OLDFIELD75}
\end{figure}

The previously described method was implemented with this data for several choices of $T$, $L$, $Q$, and $\Delta x,\Delta v,\Delta t$.  The results of a few of these runs are presented on the following pages.  We normalize the velocity domain by choosing $Q = 1$ so that velocities are only allowed in the interval $[-1,1]$.  In each of the runs, it is important to properly balance the choices of $T$ and $L$. If $L$ is taken too small or $T$ too large, the computations of the field at each gridpoint will be invalid after some small time $T_0 < T$. Thus, we must ensure that the domain of validity contains a region $[0,T] \times [-I,I]$ for some $I \leq L$ (see Figure \ref{DOV}). For each of the following runs, the mesh sizes are taken to be $\Delta t = \Delta x = \Delta v = 0.01$ and $L = 50$ so that the computations are valid over a reasonably long time interval.  More specifically, since maximum velocities are on average near $Q = 1$, we expect the time of validity to be approximately $T = L/Q \approx 50$.  Figure \ref{OLDFIELD50} displays the computed electric field and conserved energy for $T = 30$, while Figure \ref{OLDFIELD75} in the next section will display the same information for $T = 50$.

From Figure \ref{OLDFIELD50}, we can see that the computed electric
field is oscillatory, but with decreasing amplitude and mean for $0 \leq t \leq 30$.  In addition, after considering values of $t$
near which the peaks of these oscillations occur, we may conclude
this amplitude decreases around a rate of $t^{-1}$ for $t
\in [0,30]$.  This result is demonstrated in Table \ref{TABLE2}.
Specifically, the last row of the table leads us to believe that for
$t\in [0,30]$,
$$ \sup_l \vert E(t, x_l) \vert \approx 1.1 \times 10^{-4} \ t^{-1}.$$

\begin{table}[t]
\centering
\begin{tabular}{|c|c|c|c|c|}
\hline
\ & t = 15.2 & t = 19.0 & t = 23.2 & t = 26.0 \\
\hline
 $\displaystyle \sup_l \vert E(t,x_l) \vert$ & $ 6.5 \times 10^{-6}$ & $5.7 \times 10^{-6}$ & $5.3 \times 10^{-6}$ & $5.0 \times 10^{-6}$ \\
\hline
 $\displaystyle \sup_l \vert E(t,x_l) \vert \cdot t $ & $9.9 \times 10^{-5}$ & $1.1 \times 10^{-4}$  & $1.2 \times 10^{-4}$ & $1.3 \times 10^{-4}$ \\
\hline
\end{tabular}
\caption[Time rate of decay of computed field for 1-D particle
method]{Time rate of decay for the computed electric field } \label{TABLE2}
\end{table}

To the right of the computed electric field, Figure \ref{OLDFIELD50} shows the net energy in the computational domain for times $0 \leq t \leq 30$.  The deviation of the computed
energy from its initially computed value is typically seen as a good
measure of the accuracy of the method over time.  In this case, the
values of the energy range between $1.747 \times 10^{-2}$, initially,
and $1.665 \times 10^{-2}$, near time $T = 30$.  The relative change is
calculated as
$$ \mathrm{\% change} = \frac{1.755 \times 10^{-2} - 1.657 \times 10^{-2}}{1.657 \times 10^{-2}} \approx 4.92\%.$$
Thus, we expect our computations to be within $4\%$ and $5 \%$ of
their actual value.

\section*{Breakdown Outside the Domain of Validity}
\noindent While the particle method does yield accurate and efficient approximations of (\ref{VP}), limitations do exist within the formulation and implementation.  Of course, a limitation of any approximation of a particle distribution without compact support is the truncation of the spatial and velocity domains.  One can only simulate domains of finite length, and the previously described particle method will determine the true effects observed, but only within the domain of validity.  However, the size of the domain of validity is another limitation of the method.  In Figure \ref{OLDFIELD75}, we can see that after time $t = 30$, the behavior of the calculated field changes slightly.  The decay to zero is impeded and the mean of the oscillation begins to increase.  In addition, the figure demonstrates that the energy begins to decrease further away from its initial value after time $30$, as well. Within the simulation, particle velocities have increased, causing the expected domain of validity to shrink. Thus, the computations are only valid up to a time which is strictly less than the stopping time, after which other particles beginning outside of the truncated spatial domain will influence field behavior within.  This displays the limitation of the particle method, graphically represented in the previous section by Figure \ref{DOV}.  Computations will inevitably fail to be accurate on any portion of the spatial domain after some time.  The only remedy within the context of a particle method is to enlarge the initial spatial domain, resulting in more expensive computations.

\pagebreak

\end{document}